\documentclass[11pt]{article}
\usepackage{amsmath,amssymb,theorem}
\usepackage{graphicx,subfigure,enumerate}
\usepackage[ruled, lined, linesnumbered]{algorithm2e}
\usepackage{psfrag,placeins,wrapfig,booktabs,chngcntr}
\usepackage[affil-it]{authblk}
\counterwithin{figure}{section}
\numberwithin{figure}{section}
\counterwithin{table}{section}
\usepackage{hyperref}
\hypersetup{
  colorlinks,
  citecolor=black,
  filecolor=black,
  linkcolor=black,
  urlcolor=black
}

\newtheorem{thm}{Theorem}[section]
\newtheorem{conj}[thm]{Conjecture}
\newtheorem{lem}[thm]{Lemma}%[section]
\theorembodyfont{\rmfamily}
\newtheorem{prob}[thm]{Problem}%[section]
\def\pf{\bigskip\noindent {\bf Proof.}~~}

\def\dfn#1{{\sl #1}}
\def\es{\emptyset}
\def\less{\setminus}

\def\pf{\bigskip\noindent {\bf{Proof.}}~~}

\def\qed{ \hfill $\square$}

\newcounter{counter}
\def\eps{\varepsilon}

  \allowdisplaybreaks[4]

\begin{document}
\title{On the size of special class 1 graphs and     $(P_3; k)$-co-critical graphs}
\author{Gang Chen$^{1,}$\thanks{Supported by  the National   Science  Foundation of China under Grant No.  12061056.}\,\,, Zhengke Miao$^{2,}$\thanks{This work was done in part while   the second author visited the University of Central Florida in May 2019. Supported by the National   Science  Foundation of China under Grant No. 11971205.}\,\,, Zi-Xia Song$^{3,}$\thanks{Supported by the National   Science  Foundation under Grant No. DMS-1854903. E-mail address:  %chen_g@nxu.edu.cn; zkmiao@jsnu.edu.cn; 
Zixia.Song@ucf.edu.% zhangjingmei@hsas.edu.cn
} \,and Jingmei Zhang$^4$ 
}
 
   \affil{
  { \small {$^1$School of Mathematics and Statistics, 
Ningxia University, Yinchuan, Ningxia 750021, China}}\\
   { \small {$^2$Research Institute of Mathematical Science and  Department of Mathematics and Statistics, \\
Jiangsu Normal University, Xuzhou, Jiangsu 221116, China}}\\
  { \small {$^3$Department  of Mathematics, University of Central Florida, Orlando, FL 32816, USA}} \\
  { \small {$^4$ The High School affiliated to the Southern University of Science and Technology, \\
  Shenzhen, Guangdong  518133, China   }} \\
    }

\date{}
\maketitle

\begin{abstract}
A well-known theorem of Vizing  states that if $G$ is a simple graph with
maximum degree $\Delta$, then the   chromatic index  $\chi'(G)$  of $G$ is $\Delta$  or $\Delta+1$.  A
graph $G$ is    class 1 if $\chi'(G)=\Delta$, and     class 2 if $\chi'(G)=\Delta+1$; $G$ is 
 $\Delta$-critical if it is connected, class 2  and $\chi'(G-e)<\chi'(G)$ for  every  $e\in E(G)$. A long-standing  conjecture of Vizing from 1968 states that every $\Delta$-critical graph on $n$ vertices has at least  $(n(\Delta-1)+ 3)/2$ edges.  We initiate the study of determining the minimum number of edges of class 1 graphs $G$, in addition,  $\chi'(G+e)=\chi'(G)+1$ for  every  $e\in E(\overline{G})$. Such graphs have intimate relation to  $(P_3; k)$-co-critical graphs, where   a non-complete graph $G$ is $(P_3; k)$-co-critical if  there exists a $k$-coloring of $E(G)$ such that $G$ does not contain a monochromatic copy of $P_3$  but   every $k$-coloring of   $E(G+e)$ contains  a monochromatic copy of $P_3$ for every $e\in E(\overline{G})$. We  use the bound on the size of the aforementioned  class 1 graphs  to study the  minimum number of edges over all  $(P_3; k)$-co-critical graphs. We prove that if $G$ is  a  $(P_3; k)$-co-critical graph on $n\ge k+2$ vertices,    then \[e(G)\ge {k \over 2}\left(n- \left\lceil {k \over 2} \right\rceil - \varepsilon\right) + {\lceil k/2 \rceil+\varepsilon \choose 2},\]  
   where $\varepsilon$ is the remainder of  $n-\lceil k/2 \rceil $ when divided by $2$. This bound is best possible for   all $k \ge 1$ and $n \ge \left\lceil {3k /2} \right\rceil +2$.  
\end{abstract}



\baselineskip 16pt

\section{Introduction}
All graphs considered in this paper are finite, and without loops or multiple edges. For a graph $G$, we   use $V(G)$ to denote the vertex set, $E(G)$ the edge set,  $e(G)$ the number of edges,  $N(x)$   the neighborhood of vertex $x$ in $G$, $\delta(G)$ the minimum degree, $\Delta(G)$ the maximum degree,  and $\overline{G}$ the complement of $G$.
For $A, B \subseteq V(G)$, we denote by $B \less A$ the set $B - A$, 
 $e(A, B)$ the number of edges between $A$ and $B$ in $G$,  and $G \less A$ the subgraph of $G$ induced on $V(G) \less A$, respectively.
If $A = \{a\}$, we simply write $B \less a$, $e(a, B)$, and $G \less a$, respectively.  For  every $e\in E(G)$ and      $e'\in  \overline{G} $, we use $G-e$ and $G+e'$ to denote the graph obtained from $G$ by deleting the edge $e$ and the graph obtained from $G$ by   adding the new edge $e'$, respectively.
   For an integer $t\ge1$ and a graph $H$, we define $tH$ to be the union of $t$ disjoint copies of $H$.  We use $K_n$, $P_n$     and $T_n$ to denote the complete graph, a path  and a tree on  $n$ vertices, respectively.  For any positive integer $k$, we write  $[k]$ for the set $\{1,2, \ldots, k\}$. We use the convention   ``$A:=$'' to mean that $A$ is defined to be the right-hand side of the relation.  \medskip

Given an integer $k \ge 1$ and   graphs $G$, ${H}_1, \dots, {H}_k$, we write \dfn{$G \rightarrow ({H}_1, \dots, {H}_k)$} if every $k$-coloring of $E(G)$ contains a monochromatic  ${H}_i$ in color $i$ for some $i\in [k]$.
The classical \dfn{Ramsey number}  $r({H}_1, \dots, {H}_k)$  is the minimum positive integer $n$ such that $K_n \rightarrow ({H}_1, \dots, {H}_k)$.  Following Ne\v{s}et\v{r}il~\cite{Nesetril1986}, and Galluccio, Simonovits and Simonyi~\cite{Galluccio1992}, a non-complete graph $G$   is \emph{$(H_1, \ldots,  H_k)$-co-critical} if $G  \nrightarrow ({H}_1,  \ldots, {H}_k)$,   but  $G+e\rightarrow ({H}_1,  \ldots, {H}_k)$ for every edge $e$ in $\overline{G}$.    We simply write $r({H}_1, \dots, {H}_k)$  by $r(H;k)$ and say  $G$ is $(H; k)$-co-critical  when  $H_1=\cdots=H_k=H$. \medskip

 The notation of co-critical graphs was initiated by Ne\v{s}et\v{r}il~\cite{Nesetril1986} in 1986. It is simple to check that   $K_6^-$ is $(K_3, K_3)$-co-critical, where   $K_6^-$ denotes  the graph obtained from $K_6$ by deleting exactly one edge. 
It is worth noting that every   $(H_1, \ldots,  H_k)$-co-critical graph has at least $r({H}_1, \dots, {H}_k)$   vertices.      Hanson and Toft~\cite{Hanson1987} in 1987 also studied  the minimum and maximum number  of edges over all $(H_1, \ldots,  H_k)$-co-critical graphs on $n$ vertices when $H_1, \ldots,  H_k$ are complete graphs, under the name of \dfn{strongly $(|H_1|, \ldots, |H_r|)$-saturated} graphs. Recently, this topic has  been studied under the name of \dfn{$\mathcal{R}_{\min}(H_1, \dots, H_k)$-saturated} graphs \cite{Chen2011, Ferrara2014,  RolekSong, SongZhang}. We refer the reader to   a recent paper  by  the last two authors \cite{SongZhang} for
further background on $(H_1, \ldots,  H_k)$-co-critical graphs.  Hanson and Toft~\cite{Hanson1987}   made the following conjecture from 1987 on   $(K_{t_1}, \dots, K_{t_k})$-co-critical graphs.

\begin{conj}[Hanson and Toft~\cite{Hanson1987}]\label{HTC}  Let   $r = r(K_{t_1}, \dots, K_{t_k})$. Then  
 every   $(K_{t_1}, \dots, K_{t_k})$-co-critical graph   on $n$ vertices has at least    
$$   (r- 2)(n - r + 2) + \binom{r - 2}{2}$$
edges. This bound is best possible for every $n$. 
\end{conj}
\medskip

It was shown in \cite{Chen2011} that every   $(K_3,  K_3)$-co-critical graph on $n\ge 56$ vertices has at least $4n-10$ edges, thereby verifying the first nontrivial case of Conjecture~\ref{HTC}. At this time, however, it seems that a complete resolution of Conjecture~\ref{HTC} remains elusive.  
 Some structural properties of $(K_3,  K_4)$-co-critical graphs are given in \cite{K3K4}.    Motivated by Conjecture~\ref{HTC}, Rolek and the third author~\cite{RolekSong} recently initiated the study of the minimum number of possible edges over all $ (K_t, \mathcal{T}_k)$-co-critical graphs, where  $\mathcal{T}_k$ denotes  the family of all trees on $k$ vertices, and 
for all $t, k \ge 3$, we write $G\rightarrow (K_t, \mathcal{T}_k)$  if  for every $2$-coloring   $\tau: E(G) \to \{\text{red, blue} \}$, $G$ has either a red $K_t$ or a blue tree $T_k\in \mathcal{T}_k$; 
a non-complete graph $G$ is \dfn{$(K_t, \mathcal{T}_k)$-co-critical} if  
$G\nrightarrow (K_t, \mathcal{T}_k)$, but $G+e\rightarrow (K_t, \mathcal{T}_k)$ for all $e$ in $\overline{G}$.    
The following results have been   obtained on the size of $ (K_t, \mathcal{T}_k)$-co-critical graphs.     

\begin{thm}[Rolek and Song~\cite{RolekSong}]\label{K3Tk}  Let $n, k\in \mathbb{N}$.  
\begin{enumerate}[(i)]
 \item  Every  $(K_3, \mathcal{T}_4)$-co-critical graph on $n\ge 18$ vertices has  at least $\left\lfloor 5n/2\right\rfloor$ edges. This bound is sharp for every $n\ge18$.  
 \item   For all  $k \ge 5$,  if $G$ is $(K_3, \mathcal{T}_k)$-co-critical on $n\ge  2k + (\lceil k/2 \rceil +1) \lceil k/2 \rceil -2$ vertices, then $$ e(G) \ge \left(\frac{3}{2}+\frac{1}{2}\left\lceil \frac{k}{2} \right\rceil\right)n-c(k),$$ 
 where $c(k)=\left(\frac{1}{2} \left\lceil \frac{k}{2} \right\rceil + \frac{3}{2} \right) k -2$. This bound is asymptotically best possible.
 \end{enumerate}
\end{thm}

\begin{thm}[Song and Zhang~\cite{SongZhang}]\label{SongZhang}   
  Let  $ t, k\in \mathbb{N}$ with $t \ge 4$ and  $k\ge\max\{6, t\}$.  There exists a  constant  $\ell(t, k)$  such that, for  all $n\in \mathbb{N}$ with    $n \ge (t-1)(k-1)+1$,    if $G$  is  a $(K_t, \mathcal{T}_k)$-co-critical graph  on $n$ vertices, then     $$   e(G)\ge   \left(\frac{4t-9}{2}+\frac{1}{2}\left\lceil \frac{k}{2} \right\rceil\right)n-\ell(t, k).$$ 
This bound is asymptotically best possible    when $t\in\{4,5\}$ and $k\ge6$.  
 \end{thm}

 The methods developed in \cite{RolekSong, SongZhang}  may shed some light on attacking Conjecture~\ref{HTC}.  Inspired  by Conjecture~\ref{HTC},  Ferrara,    Kim and  Yeager~\cite{Ferrara2014} proposed the following problem.
 
 \begin{prob}[Ferrara,    Kim and  Yeager~\cite{Ferrara2014}]\label{Problem} Let $H_1,   \ldots,  H_k$ be graphs, each with   at least one edge. Determine  the minimum number of edges of $(H_1,   \ldots,  H_k)$-co-critical graphs.
 \end{prob}
 
  In the same paper they settled Problem~\ref{Problem} when each $H_i$ is a matching of $m_i$ edges. 
  
  \begin{thm}[Ferrara,    Kim and  Yeager~\cite{Ferrara2014}]\label{matching}
Let $m_1,   \ldots, m_k$ be positive integers. Then every $(m_1K_2,   \ldots,  m_kK_2)$-co-critical graph on $n>3(m_1+\cdots +m_k-k)$ vertices has at least $3(m_1+\cdots +m_k)$ edges. This bound is best possible for all $n>3(m_1+\cdots +m_k-k)$. 
\end{thm}
 
 Theorem~\ref{matching} yields the very first result  on Problem~\ref{Problem} for multicolor  $k$.   In this paper, we    continue to study Problem~\ref{Problem} by determing the minimum number of possible edges over all $ (P_3; k)$-co-critical graphs for multicolor $k$. It was shown  in \cite{paths} that    
$r(P_3; k)=2\lceil k/2\rceil +1 $ for all $k \ge 1$. Thus every $ (P_3; k)$-co-critical graph has at least $2\lceil k/2\rceil +1$ vertices. Before we state our main results, we need to introduce more notation. For a positive integer $t$, a graph $H$ is \dfn{properly $t$-edge-colorable} if there exists a mapping  $\varphi: E(H) \rightarrow [t]$ such that $\varphi(e)\ne \varphi(f)$ for every pair of adjacent edges  $e$ and $f$   in $H$. The \dfn{chromatic index} $\chi'(H)$ of $H$ is the least   $t$ such that $H$ is  properly $t$-edge-colorable. 
The following is a  well-known result of Vizing~\cite{Vizing68}. 

\begin{thm}[Vizing~\cite{Vizing68}]\label{Vizing}  Every graph $H$ satisfies 
 
 \[ \Delta(H)\le\chi'(H)\le \Delta(H)+1.\]
 
 \end{thm}

A
graph $H$ is   class 1 if $\chi'(H)=\Delta(H)$, and    class 2 if $\chi'(H)=\Delta(H)+1$. A connected graph $H$ is 
 $\Delta$-critical if $\Delta(H)=\Delta$, class 2  and $\chi'(H-e)<\chi'(H)$ for  every  $e\in E(H)$.  Let $G$ be a  $ (P_3; k)$-co-critical graph and let   $\tau: E(G) \rightarrow [k]$ be a $k$-coloring of $E(G)$. We say that  $\tau$ is a \dfn{critical coloring} of $G$ if 
$G $ has no monochromatic copy of $P_3$ under $\tau$. Since $G$ is  $ (P_3; k)$-co-critical, we see that  $G$ admits a critical coloring but $G+e$ has no critical coloring for every $e\in E(\overline{G})$. Let $\tau: E(G) \rightarrow [k]$ be  a critical coloring of $ G$. Then $G$ has no monochromatic copy of $P_3$ under the coloring $\tau$. It follows that:  adjacent edges of $G$ are colored differently under $\tau$, and so  $\tau$ is a proper $k$-edge-coloring of $G$;    $G+e$  is not properly  $k$-edge-colorable for every $e\in E(\overline{G})$. This motivates us to study the minimum number of edges of class 1 graphs $H$ such that $\chi'(H+e)=\chi(H)+1$ for every $e\in E(\overline{H})$. Note that $H+e$ may not be of class 2 for every $e\in E(\overline{H})$). It is also worth noting that a long-standing conjecture of  Vizing~\cite{Vizing68}  from 1968 states that every $\Delta$-critical graph on $n$ vertices has at least  $(n(\Delta-1)+ 3)/2$ edges. We prove the following main result. 
\medskip

 \begin{thm}\label{Song}
Let $H$  be a class 1 graph with $n$ vertices and maximum degree $\Delta$. If   $\chi'(H+e)=\Delta +1$ for every $e\in E(\overline{H})$, then 
 $$e(H)\ge {\Delta\over 2}\left(n- \left\lceil {\Delta \over 2} \right\rceil - \varepsilon\right) + {\lceil \Delta/2 \rceil+\varepsilon \choose 2},$$
   where $\varepsilon$ is the remainder of $n-\lceil \Delta/2 \rceil$ when divided by $2$. This bound is best possible  for all $\Delta\ge1$ and $n \ge \lceil 3\Delta/2 \rceil $. 
\end{thm}
 
The study of the minimum number of edges of class 1 graphs $H$ such that $\chi'(H+e)=\chi'(H)+1$ for every edge $e\in E(\overline{H})$  may be of independent interest.   We prove Theorem~\ref{Song} in Section~\ref{class1}. We then apply Theorem~\ref{Song} to establish the minimum number of edges over all  $ (P_3; k)$-co-critical graphs in Section~\ref{co-critical}. We prove the following. 

 \begin{thm}\label{Lower}
For all $k \ge 1$ and $n \ge2\lceil k/2\rceil +1$, every  $(P_3;k)$-co-critical graph $G$ on $n$ vertices satisfies  
 $$e(G)\ge {k \over 2}\left(n- \left\lceil {k \over 2} \right\rceil - \varepsilon\right) + {\lceil k/2 \rceil+\varepsilon \choose 2},$$
   where $\varepsilon$ is the remainder of  $n-\lceil k/2 \rceil $ when divided by $2$. This bound is best possible for   all $k \ge 1$ and $n \ge \left\lceil {3k /2} \right\rceil +2$. 
\end{thm}

We end this section with  the well-known Vizing's Adjacency Lemma~\cite{Vizing65} that will be needed in the proofs of our main results.   %An edge $e$ of   a class 2 graph $H$ is  a \dfn{critical edge} if $\chi'(H- e)<\chi'(H)$. 

 \begin{lem}[Vizing's Adjacency Lemma~\cite{Vizing65}]\label{VAL}    
Let $H$ be a class 2 graph and let $xy\in E(H)$ be a critical edge, that is,  $\chi'(H- xy)<\chi'(H)$. Then
  $x$ is adjacent to at least $\Delta(H)+1-d(y)$ vertices of degree $\Delta(H)$ in $H- xy$.  In particular,   $H$ has at least  three  vertices of degree $\Delta(H)$. 
 \end{lem}

%%%%%%%%%%

\section{Proof of Theorem~\ref{Song}}\label{class1}

   Let $H$  be a class 1 graph with $n$ vertices and   $\Delta(H)=\Delta $ such that     $\chi'(H+e)=\Delta +1$ for every $e\in E(\overline{H})$.  
  It is simple to check that the statement is true when $\Delta\le 1$. We may assume that $\Delta\ge2$.  Suppose 
\begin{align*}
e(H)< {\Delta \over 2}\left(n- \left\lceil {\Delta \over 2} \right\rceil - \varepsilon\right) + {\lceil \Delta/2 \rceil+\varepsilon \choose 2} = {\Delta n \over 2}  -{1\over 2} \left(\left\lceil {\Delta \over 2} \right\rceil +\eps\right)\left( \left\lfloor {\Delta\over2}\right\rfloor-\eps+1\right).\end{align*}
Then 
 \begin{align*}2e(H) &< \Delta n-\left(\left\lceil {\Delta  \over 2} \right\rceil +\eps\right)\left( \left\lfloor {\Delta \over2}\right\rfloor-\eps+1\right) 
 =\begin{cases}  \Delta n- \frac{(\Delta +1)^2}4+\frac14  &\mathrm{if} \,\, \Delta \text{ is even}\\
\Delta n- \frac{(\Delta +1)^2}4   &\mathrm{if} \,\, \Delta \text{ is odd and } \eps=0\\
\Delta n- \frac{(\Delta +1)^2}4+1  &\mathrm{if} \,\, \Delta \text{ is odd and } \eps=1
\end{cases} \label{edges}\tag{$*$} 
\end{align*}
Let $\delta:=\delta(H)$. Then $\delta  \le \Delta-1$  by (\ref{edges}).  Let   $\varphi: E(H) \rightarrow [\Delta]$ be a proper $\Delta$-edge-coloring of  $H$. For each $i \in [\Delta]$, let $E_i:=\{uv \in E(H)\mid  \varphi(uv)=i\}$.  For   $i, j \in [\Delta]$ with $i \ne j$, let  $H_{i,j}$ be  the spanning  subgraph of $H$ with edge set $E_i \cup E_j$.  Then each component of $H_{i,j}$ is either an even cycle or a path. 
For each  $u \in V(H)$, let $\varphi(u):=\{\varphi(uv)\mid  v \in N(u)\}$ and $\bar{\varphi}(u)=[\Delta] \less \varphi(u)$.
If $\alpha \in \varphi(u)$ and $\beta \in \bar{\varphi}(u)$, then the component of $H_{\alpha, \beta}$ containing $u $ is a path,    and called the \dfn{$(\alpha, \beta)$-Kempe chain} (or simply \dfn{$(\alpha, \beta)$-chain}) starting at $u$.  
We next prove several claims.\bigskip

\setcounter{counter}{0}

\noindent {\bf Claim\refstepcounter{counter}\label{delta}  \arabic{counter}.}  
$\delta\ge1$.

\pf Suppose $H$ has   an isolated vertex, say $u$. Then there is a vertex $v \in V(H) $ with $v\ne u$ such that $d(v)\le \Delta-1$, else $e(H)=\Delta(n-1)/2 $,  contrary to (\ref{edges}). Note that $uv\notin E(H)$. But then we obtain a proper $\Delta$-edge-coloring of $H+uv$  from $\varphi$ by coloring the edge $uv$ with a color in $\bar{\varphi}(v)$, which contradicts to the assumption that $\chi'(H+e)=\Delta +1$ for every $e\in E(\overline{H})$. \qed\bigskip

Let    $A:=\{v \in V(H)\mid d(v) \le \Delta-1\}$ and $B:=V(H)\less A$. 
Then $A\ne \emptyset $ and  $B\ne \emptyset$.  Furthermore, for each $v\in A$, we have $\varphi(v)\ne \emptyset$ by Claim~\ref{delta},  and $\bar{\varphi}(v)\ne \emptyset$ by the definition of $A$. 
 \bigskip

\noindent {\bf Claim\refstepcounter{counter}\label{kempechain}  \arabic{counter}.}  
 Let  $u, v \in A$ with $uv \not\in E(H)$. Then 
\begin{enumerate}[(i)]
\item $\Delta(H+uv)=\Delta$ and $H+uv$ is of class 2. Furthermore,    $uv$ is a critical edge in $H+uv$.

\item   for every  $\alpha \in \varphi(u) \less \varphi(v)$ and $\beta \in \varphi(v) \less \varphi(u)$,   the  $(\alpha, \beta)$-chain starting at $u$  terminates at  $v$.  
 \item  $\bar{\varphi}(u)\subseteq  \varphi (v)$ and  $\bar{\varphi}(v)\subseteq  \varphi (u)$.   Moreover,  $d(u)+d(v)\ge \Delta$.
  \end{enumerate}

\pf To prove (i),   let $G:=H+uv$. Since $u, v\in A$, we see that  $d_{G}(u)\le \Delta$ and $d_{G}(v)\le \Delta$. Note that $\chi'(G)=\Delta +1$ by assumption and   $\Delta(G)=\Delta$.  Thus  $G$ is of class 2  and $uv$ is a critical edge in $G$ because  $\chi'(G-uv)=\chi'(H)=\Delta$.   \bigskip

 To prove (ii), 
let $\alpha \in \varphi(u) \less \varphi(v)$ and $\beta \in \varphi(v) \less \varphi(u)$.       Let $J$ be the component of $H_{\alpha, \beta}$ that contains the vertex  $u$. Then $J$ is a path. Suppose $v \not\in V(J)$.  Let $\xi$ be a proper $\Delta$-edge-coloring of $H$ obtained from $\varphi$ by interchanging the colors $\alpha$ and $\beta$ in $J$. Note that $\alpha \in \bar{\xi}(u) \cap \bar{\xi}(v)$. Then we obtain a proper $\Delta$-edge-coloring of  $H+uv$  from $\xi$ by coloring the edge $uv$ with color $\alpha$, contrary to Claim~\ref{kempechain}(i). This proves that  $v\in V(J)$.   Then $J$ is a path with ends $u, v$. Therefore,    $J$ is the desired $(\alpha, \beta)$-chain.  \bigskip

 It remains to prove (iii).   Note that    $ \bar{\varphi}(u) \cap \bar{\varphi}(v)=\es$, else we obtain a proper $\Delta$-edge-coloring of  $H+uv$ from $\varphi$ by coloring the edge $uv$ with a color in   $  \bar{\varphi}(u) \cap \bar{\varphi}(v)$, contrary to Claim~\ref{kempechain}(i). Thus $\bar{\varphi}(u)\subseteq  \varphi (v)$ and  $\bar{\varphi}(v)\subseteq  \varphi (u)$. Note that  $\Delta\ge | \bar{\varphi}(u)|+ |\bar{\varphi}(v)|=(\Delta-d(u))+(\Delta-d(v))$, and so $d(u)+d(v)\ge \Delta$.  
 \qed\bigskip 
 
 Let $x\in A$ with $d(x)=\delta$.   Let  $N[x]:=N(x)\cup\{x\}$,   $\ell:=|A\less N[x]|\ge0$, and $  A\less N[x]:=\{y_1, y_2, \ldots, y_\ell\}$ when $\ell\ge1$.  \bigskip

\noindent {\bf Claim\refstepcounter{counter}\label{Alpha}  \arabic{counter}.}   If $\ell\ge1$, then for every $\alpha\in\varphi(x)$,  we have   $\alpha\notin  \bar{\varphi}(y_i) \cap \bar{\varphi}(y_j)$ for every pair of  distinct $i, j\in[\ell]$. Thus $\bar{\varphi}(y_1),    \bar{\varphi}(y_2), \ldots, \bar{\varphi}(y_\ell)$ are pairwise disjoint  and $ \sum_{i=1}^\ell |\bar{\varphi}(y_i) |\le d(x)$.

\pf  Suppose  there exists some $\alpha\in\varphi(x)$ such that  $\alpha\in  \bar{\varphi}(y_i) \cap \bar{\varphi}(y_j)$ for some    $i, j\in[\ell]$ with $i\ne j$.    Let $\beta\in \bar{\varphi}(x)$. By Claim~\ref{kempechain}(iii), $\beta\in \varphi(y_i) \cap \varphi(y_j)$.   By Claim~\ref{kempechain}(ii), let $P$ be  the  $(\alpha, \beta)$-chain starting at $x$ and ending at $y_i$. Note that  $y_j\notin V(P)$ because $\alpha\in  \bar{\varphi}(y_j)$.  But then we   obtain a proper $\Delta$-edge-coloring of $H+xy_j$ from $\varphi$ by first   interchanging the colors $\alpha$ and $\beta$  on  $P$ and then   coloring the edge $xy_j$ with   color   $\alpha$, contrary to Claim~\ref{kempechain}(i). This proves that    $\bar{\varphi}(y_1),    \bar{\varphi}(y_2), \ldots, \bar{\varphi}(y_\ell)$ are pairwise disjoint. By Claim~\ref{kempechain}(iii), $\bar{\varphi}(y_i) \subseteq  \varphi(x)$ for each $i\in[\ell]$.    Thus  $ \sum_{i=1}^\ell |\bar{\varphi}(y_i) |\le d(x)$.\qed\bigskip

\noindent {\bf Claim\refstepcounter{counter}\label{e(x,B)}  \arabic{counter}.}  
If $\ell\ge1$, then     $e(x, B)\ge  1$.
 
 \pf Suppose $\ell\ge1$. Then   $y_1\in A\less N[x]$ and $xy_1\notin E(H)$.   By Claim~\ref{kempechain}(i), $xy_1$ is a critical edge of $H+xy_1$. By Lemma~\ref{VAL}, $x$ is adjacent to at least $\Delta+1-d_{H+xy_1}(y_1)=\Delta-d(y_1) $ vertices of degree $\Delta$ in $H$.  Hence  $e(x, B)\ge \Delta-d(y_1)\ge1$.\qed\bigskip
 
 We complete the proof by considering the following two cases on $\ell$. 
%\noindent {\bf Claim\refstepcounter{counter}\label{E(H)}  \arabic{counter}.}  $2e(H)\ge   \begin{cases}  \Delta n- \frac{(\Delta +1)^2}4+\frac14  &\mathrm{if} \,\, \Delta \text{ is even}\\
%\Delta n- \frac{(\Delta +1)^2}4  &\mathrm{if} \,\, \Delta \text{ is odd and } \eps=0\\
%\Delta n- \frac{(\Delta +1)^2}4+1  &\mathrm{if} \,\, \Delta \text{ is odd and } \eps=1
%\end{cases}   $. 
 We first consider the case   $\ell=0$. Then $A\subseteq N[x]$ and so $|A|\le |N[x]|\le \delta+1$. Thus 
\begin{align*}
2e(H) \ge \delta\cdot |A|+\Delta(n-|A|) 
&=\Delta n-(\Delta-\delta)|A|\\
&\ge \Delta n-(\Delta-\delta)(\delta+1)\\
&=\Delta n+\delta^2-(\Delta -1)\delta-\Delta \\
&=\Delta n+\left(\delta-\frac{\Delta -1}2\right)^2 -\frac{(\Delta +1)^2}4.  \end{align*} 
Since $2e(H)$ is even, by considering the parity of $n$, $\Delta$ and $\lceil(\Delta+1)/2\rceil$, it follows that 
\[2e(H) \ge  \begin{cases}  \Delta n- \frac{(\Delta +1)^2}4+\frac14  &\mathrm{if} \,\, \Delta \text{ is even}\\
\Delta n- \frac{(\Delta +1)^2}4  &\mathrm{if} \,\, \Delta \text{ is odd and } \eps=0\\
\Delta n- \frac{(\Delta +1)^2}4+1  &\mathrm{if} \,\, \Delta \text{ is odd and } \eps=1,\end{cases} \]  
  which contradicts   (\ref{edges}).  \medskip

It remains to  consider the case  $\ell\ge1$.  By Claim~\ref{e(x,B)}, $e(x, B)\ge1$ and so   $|A\cap N[x]|\le \delta$.       By Claim~\ref{Alpha},   $\bar{\varphi}(y_1),    \bar{\varphi}(y_2), \ldots, \bar{\varphi}(y_\ell)$ are pairwise disjoint and  $ \sum_{i=1}^\ell|\bar{\varphi}(y_i)  |\le d(x)=\delta$. It follows that 
\begin{align*}
\sum_{i=1}^\ell d(y_i) &=\sum_{i=1}^\ell |\varphi(y_i)| 
 =\sum_{i=1}^\ell (\Delta -|\bar{\varphi}(y_i)|)
 =\Delta \ell -\sum_{i=1}^\ell|\bar{\varphi}(y_i)  | 
 \ge \Delta \ell  -\delta. 
\end{align*}
 Note that $|A|= |A\cap N[x]|+\ell $.  Then  
\begin{align*}
2e(H) \ge \delta\cdot |A\cap N[x]|+\sum_{i=1}^\ell d(y_i) +\Delta (n-|A|)  
&\ge\delta\cdot |A\cap N[x]|+(\Delta \ell-\delta)+\Delta (n- |A\cap N[x]|-\ell)\\
&=\Delta n-\delta-(\Delta -\delta) |A\cap N[x]|\\
&\ge \Delta n -\delta-(\Delta -\delta)   \delta \\
&=\Delta n+\delta^2-(\Delta +1)\delta\\
&=\Delta n+\left(\delta-\frac{\Delta +1}2\right)^2 -\frac{(\Delta +1)^2}4.  
\end{align*} 
It follows that 
\[ 2e(H) \ge \begin{cases}  \Delta n- \frac{(\Delta +1)^2}4+\frac14  &\mathrm{if} \,\, \Delta \text{ is even}\\
\Delta n- \frac{(\Delta +1)^2}4  &\mathrm{if} \,\, \Delta \text{ is odd and } \eps=0\\
\Delta n- \frac{(\Delta +1)^2}4+1  &\mathrm{if} \,\, \Delta \text{ is odd and } \eps=1, 
\end{cases} \]     
 which   contradicts (\ref{edges}).  \medskip

 This completes the proof of Theorem~\ref{Song}. The sharpness of the bound in Theorem~\ref{Song} follows from    Lemma~\ref{upper} in Section~\ref{co-critical}. \qed\bigskip

\section{Proof of Theorem~\ref{Lower}}\label{co-critical}
We begin this section with a lemma which shows that  the  bounds  for both Theorem~\ref{Song} and Theorem~\ref{Lower} are best possible.

\begin{lem}\label{upper}
For all $k \ge 1$ and $n \ge  \left\lceil {3k / 2} \right\rceil  +1+\varepsilon$, there exists a $(P_3; k)$-co-critical graph $G$ on $n$ vertices such that 
 $$e(G)={k \over 2}\left(n- \left\lceil {k \over 2} \right\rceil - \varepsilon\right) + {\lceil k/2 \rceil+\varepsilon \choose 2},$$
%\[
%e(G) = \begin{cases}
%(4kn-k^2-2k+4\varepsilon-1)/8 & \text{ if } k  \text{ is odd},\\
%(4kn-k^2-2k)/8 &  \text{ if }  k \text{ is even}.
% \end{cases}
%\]
  where $\varepsilon$ is the remainder of  $n-\lceil k/2 \rceil$ when divided by $2$.  \end{lem}

\begin{pf} Let $k, n, \eps$ be as given in the statement.  Then $n-\lceil k/2 \rceil - \varepsilon$ is even and $n-\lceil k/2 \rceil - \varepsilon\ge k+1$.  Let $J$ be a $k$-regular graph on $n-\lceil k/2 \rceil - \varepsilon $ vertices  with $\chi'(J)=k$. Such a graph $J$ exists because  every $k$-regular bipartite graph has chromatic index $k$.    Let  $G $ be the disjoint union of $J$ and  $ K_{\lceil k/2 \rceil+\varepsilon}$. Then $\chi'(G)=k$ and $e(G)={k \over 2}\left(n- \left\lceil {k \over 2} \right\rceil - \varepsilon\right) + {\lceil k/2 \rceil+\varepsilon \choose 2}$. Let  $e\in E(\overline{G})$. Then $  \Delta(G+e)=k+1$ and $G+e$ has at most two vertices of  degree $\Delta(G+e)$. By Lemma~\ref{VAL}, $G+e$ is not class 2 and so $\chi'(G+e)= \Delta(G+e)=k+1$. It follows that  every $k$-coloring of $E(G+e)$ contains a monochromatic copy of $P_3$. Therefore,  $G$ is $(P_3;k)$-co-critical, as desired.\qed \end{pf} \bigskip

We are now ready to prove Theorem~\ref{Lower}.  Let $k, n, \eps, G$ be as stated in Theorem~\ref{Lower}.   Let $\tau: E(G) \rightarrow [k]$ be a critical coloring of $G$. Since $G$ has no monochromatic $P_3$ under $\tau$, we see that $\tau $ is a proper $k$-edge-coloring of $G$.    By Theorem~\ref{Vizing},   $\Delta(G) \le k$.     Let  $e\in E(\overline{G})$.  Then    $G+e$  has no critical coloring because $G$ is $(P_3; k)$-co-critical. Thus  $ G+e $  is not properly $k$-edge-colorable and so $\chi'(G+e)\ge k+1$. Note that $\Delta(G+e)\le \Delta(G)+1$ and    $G+e$ has at most two vertices of degree $\Delta(G)+1$. By Lemma~\ref{VAL},  $G+e$ is not class 2 and so $\chi'(G+e)\le \Delta(G)+1\le k+1$. It follows that $\chi'(G+e)=k+1$ and $\Delta(G)=k$. Thus $\chi'(G)=\Delta(G) = k$, and $\chi'(G+e)= k+1$ for every $e\in E(\overline{G})$.  By Theorem~\ref{Song} applied to $G$, we have  
 $$e(G)\ge {k \over 2}\left(n- \left\lceil {k \over 2} \right\rceil - \varepsilon\right) + {\lceil k/2 \rceil+\varepsilon \choose 2},$$
as desired.  By  Lemma~\ref{upper}, this bound  is best possible for all $k\ge 1$ and $n\ge \left\lceil {3k / 2} \right\rceil +2$. \medskip

 This completes the proof of Theorem~\ref{Lower}. \qed\medskip

\section{Concluding remarks}

   Let $t_1, t_2, \ldots, t_k$ be positive integers.   It would be interesting to know the  minimum number of edges of $(t_1P_3, t_2P_3, \ldots, t_kP_3)$-co-critical graphs.  Theorem~\ref{Lower} gives an answer to the case when $t_1=t_2=\cdots=t_k=1$.  It would also be interesting to know if the bound in Theorem~\ref{Lower} is best possible when $n< \left\lceil {3k / 2} \right\rceil +2$.

\end{document}